\documentclass[11pt,letter]{article}

\usepackage{amssymb,epsfig,latexsym}
\usepackage{graphicx,amsmath,amssymb,latexsym,epsfig}
\usepackage{psfig,epsf}

\def\cX{{\cal X}}

\def\cX{{\cal X}}

\def\Z{\mathbb{Z}}
\def\C{\mathbb{C}}
\def\N{\mathbb{N}}
\def\R{\mathbb{R}}
\def\T{\mathbb{T}}

\def\and{\quad\mbox{and}\quad}
\def\ind{{\rm 1\hspace{-0.90ex}1}}
\def\bp{\noindent{\it Proof.}\ }
\def\ep{\hfill $\Box$}
\def\PP{{\mathbb P}}
\def\EE{{\mathbb E}}
\def\sinc{\mbox{sinc}}
\def\tr{\mbox{tr}}

\newtheorem{theorem}{\bf Theorem}
\newtheorem{proposition}[theorem]{\bf Proposition}

\newtheorem{lemma}[theorem]{\bf Lemma}

\title{Eigenvalues of Euclidean Random Matrices}
\author{Charles Bordenave \footnote{Institut de Math\'ematiques - Universit\'e de Toulouse \& CNRS (UMR 5219) -
31062 Toulouse - France. Email: bordenave@math.ups-tlse.fr }}

\begin{document}

\maketitle

\begin{abstract}

We study the spectral measure of large Euclidean random matrices. The entries of these matrices are determined by the relative position of $n$ random points in a compact set $\Omega_n$ of $\R^d$. Under various assumptions we establish the almost sure convergence of the limiting spectral measure as the number of points goes to infinity. The moments of the limiting distribution are computed, and we prove that the limit of this limiting distribution as the density of points goes to infinity has a nice expression. We apply our results to the adjacency matrix of the geometric graph.

\noindent {\bf Keywords:} random matrix, spectral measure, random geometric graphs, spatial point process, Euclidean distance matrix.

\end{abstract}

\section{Introduction}

The main research effort in the theory of random matrices concerns matrices where the coefficients are independent random variables (see Bai \cite{bai} for a survey). Few authors have studied the limiting spectral measures of other types of large matrices, in particular, Markov, Hankel and Toeplitz matrices have been studied by Bryc, Dembo and Jiang \cite{dembo2006} and Toeplitz matrices by Hammond and Miller \cite{hammond}. In this paper, we consider another class of random matrices, the Euclidean random matrices (ERM) which have been introduced by M\'ezard, Parisi and Zee \cite{mezardparisi}. An ERM is an $n \times n$ matrix, $A$,  whose entries is a function of the positions of $n$ random points in a compact set $\Omega$ of $\R^d$. In this paper, $\Omega$ will be an hypercube, the $n$ points ${\cX_n} = \{X_1,\cdots, X_n\}$, $n$ uniformly distributed points in $\Omega$ and 
\begin{equation}
\label{eq:A}
A = (F(X_i - X_j))_{ 1 \leq i \leq j \leq n},
\end{equation}
where $F$ is a measurable mapping from $\R^d$ to $\C$. We will pay attention to the spectral properties of $A$. In this paper, we will compute some limits of the spectral measure as the number of points $n$ goes to infinity. We will show how the eigenvalues of $A$ are related to the Fourier transform of the mapping $F$.

Examples of interests in branches of physics are explained in \cite{mezardparisi} and Offer and Simons \cite{offersimons}. A particularly appealing case is $F(x) = \| x\|$, the Euclidean norm. This subclass of ERM is called the random Euclidean Distance Matrices and some of their spectral properties are derived in  Vershik \cite{vershik}, Bogomolny,  Bohigas and Schmidt \cite{bogomolny}. In \cite{vershik} the author generalizes the problem considered here and consider an integral operator defined on a metric space. Koltchinskii and Gin\'e \cite{gine} have analyzed the convergence of the spectra of the matrix $(n ^{-1} h(X_i,X_j) )_{ 1 \leq i,j \leq n}$ to the spectra of the compact integral operator with symmetric kernel $h$.

Another field of application is graph theory. Indeed, if $F(x) = \ind(0 \leq \| x\| \leq r)$, then $A$ is the adjacency matrix of the proximity (or geometric) graph (refer to Penrose \cite{penrose}). More generally if $F (X) = F(-X) \in \{0,1\}$ then $A$ is the adjacency matrix of a random  graph. The spectral properties of the adjacency matrix or related matrices are of prime interest in graph theory. For example the probability of hitting times of random walks  on graphs  is governed by the spectrum of the transition matrix (for a survey on this subject, see e.g. Section 3 in Lov{\'a}sz \cite{Lovasz}). Or, in network epidemics, the time evolution of the infected population is also closely related to the spectral radius and the spectral gap of the adjacency matrix, see Draief, Ganesh and Massouli\'e \cite{draief}. For Erd{\'o}s-Renyi random graphs, some properties of the spectrum can been computed thanks to the seminal work Wigner of \cite{wigner} and F{\"u}redi and Koml{\'o}s \cite{furedikomlos}. For power law graphs and related graphs, see Chung, Lu and Vu \cite{fanchung03}, \cite{fanchung04}.

Various generalizations of (\ref{eq:A}) would be worth to consider. Some extra randomness in the model could be added, and the entry of the matrix $i,j$ could be equal to $F_{ij} ( X_i - X_j)$, where $(F_{ij})_{1 \leq i, j \leq n}$ are i.i.d. mappings independent of the point set ${\cX_n}$. Falls into this framework the adjacency matrix of a random graph where there is an edge between two points with a probability which is deterministic function of their distance, such as the small world graphs (see for example Ganesh and Draief \cite{ganesh05}).

Another generalization is the original model of M\'ezard, Parisi and Zee \cite{mezardparisi} where  the entry $i,j$ is equal to
$$F(X_i - X_j) - u \delta_{ij} \sum_k F(X_i-X_k),$$
where $\delta_{ij}$ is the Kronecker symbol and $ u \in \R$. The case $u = 1$ is of particular interest, the matrix is then a Markov matrix.

In order to obtain the adjacency matrix of more sophisticated geometric graphs, such as the Delaunay triangulation, it would be necessary to consider an entry $i,j$ which depends on the whole point set ${\cX_n}$ and not only on $X_i - X_j$.

We will consider two models in this note. In the first model, $\Omega = [-1/2,1/2]^d$ and $F$ is $1$-periodic function: if $x,y \in \R^d$ and $x-y \in \Z^d$ then $F(x) = F(y)$. Equivalently, the point set ${\cX_n} = \{X_1,\cdots, X_n\}$ could be on the unit torus  $\T^d =  \R^d \backslash \Z^d$. We choose a periodic function in order to avoid all boundary effects with the hypercube $\Omega$. The matrix $A$ is defined by (\ref{eq:A}), where $F$ is a measurable function from $\R^d$ to $\C$. The discrete Fourier transform of $F$ is defined for all $k \in \Z^d$ by $\hat F(k) = \int_{\Omega} F(x) e^{-2i\pi k.x } dx$.  Throughout the paper, we assume that a.e. and at $0$, the Fourier series of $F$  is equal to F: 
$$
F(x) = \sum_{k \in \Z^d}   \hat F(k) e^{2i\pi k.x}.
$$
A sufficient condition is $\sum_{k \in \Z^d} |\hat F(k)| < \infty$ and $F$ continuous at $0$. This Fourier transform plays an important role in the spectrum of $A$. As an example, consider $U = (U_i)_{1 \leq i \leq n}$ a vector in $\C^d$ and assume $F$ hermitian ($F(-x) = \bar F(x)$), then a.s.
\begin{eqnarray*}
U^* A U  =   \sum_{i,j} F( X_i - X_j) U_i \bar U_j
 =  \sum_{i,j} \sum_k  \hat F(k)  e^{2i\pi k.(X_i - X_j)} U_i \bar U_j
 =    \sum_k  \hat F(k)  \Bigm| \sum_{i=1}^n e^{2i\pi k.X_i} U_i \Bigm|^2.
\end{eqnarray*}
Therefore $A$ is positive if and only if for all $k \in \Z^d$, $\hat F(k) \geq 0$.

We will compute explicitly the spectral measure  of the matrix $A_n = A/n$ as $n$ tends to $\infty$,
$$
\mu_n = \sum_{i =1 } ^ n \delta_{ \lambda_i (n) / n},
$$
where $\{ \lambda_i (n) \}_{ 1 \leq i \leq n}$ is the set of eigenvalues of $A$. Notice that $\{ \lambda_i (n) / n\}_{ 1 \leq i \leq n}$ is the set of eigenvalues of $A_n$. We define the measure:
$$\mu = \sum_{k \in \Z^d} \delta_{ \hat F (k)}.$$
Since $\lim_{ \|k\| \to \infty} \hat F (k) = 0$, $\mu$ is a counting measure with an accumulation point at $0$.

\begin{theorem}
\label{th:main}
For all Borel sets $K$ with $\mu(\partial K) = 0$ and $0 \notin \bar K$, a.s.
\begin{equation}
\label{eq:weakmu}
\lim_n \mu_n(K) = \mu(K).
\end{equation}
\end{theorem}
The convergence of the spectral measure $\mu_n$ follows also from Theorem 3.1 in \cite{gine}. As an immediate corollary, we obtain the convergence of the spectral radius of $A_n$, almost surely,
$
\lim_{n \to \infty} \max_{1 \leq i \leq n} \frac {|\lambda_i (n) |}{n} =
\max_{k \in \Z^d} |\hat F(k)|.
$
For example if $F(x) = \ind ( \max_{1 \leq i \leq d} |x_i | \leq r)$ then $\hat F(k) = r^d \prod_{i=1} ^ d \sinc ( 2 \pi k_j r)$, where $\sinc (x) = \sin (x) / x$ and $k = (k_1,\cdots, k_d) \in \Z^d$. The spectral radius of $A_n$ converges a.s. to $r^d$ and the second largest eigenvalue to $r^d \sinc (2 \pi r)$ if $r$ is small enough, thus the spectral gap is equivalent to $r^{d+2} (2\pi)^ 2 / 3!$ as $r$ goes to $0$.

Our second model is more challenging. Now, ${\cX_n} = \{X_1,\cdots, X_{n}\}$ is the set of $n$ independent points uniformly distributed on the hypercube $\delta_n^{-1} \Omega = [-\delta_n^{-1}/2,\delta_n^{-1}/2]^d $ where $\delta_n$ goes to $0$. In this second model, we scale jointly the number of points and the space. We assume that for some $\gamma > 0$,
\begin{equation}
\label{eq:delta}
\lim_n \delta_n ^d n = \gamma.
\end{equation}
$\gamma$ is the asymptotic density of the point set ${\cX_n }$. Let $f$ be a measurable function from $\R^D$ to $\C$ with support included in $\Omega$, the matrix $A$ is defined by (\ref{eq:A}) (with $F$ replaced by $f$).

Considering the change of variable $x \mapsto \delta x$, the matrix $A$ is equal to the matrix $B_n$ defined by
$$B_n = (f_{\delta_n}
(X_i - X_j))_{ 1 \leq i \leq j \leq n},$$ where $f_\delta : x \mapsto f(x/\delta)$ and the point set ${\cX_n} = \{X_1,\cdots, X_{n}\}$ is a set of $n$  independent  points uniformly distributed on $\Omega$. The spectrum of $B_n$ is denoted by
$(\lambda'_1 (n), \cdots, \lambda'_ {n} (n) )$, we define the empirical
measure of its eigenvalues:
$$
\nu_{n} = \frac {1 }{n} \sum_{i=1} ^{n} \delta_{\lambda'_i (n)},
$$
We will prove the following:
\begin{theorem}
\label{th:scmain} For all $\gamma > 0$, there exists a measure
$\nu_\gamma$ such that for the topology of the weak convergence,
a.s.:
$$
\lim_{n \to \infty} \nu_{n} = \nu_\gamma.
$$
Moreover $\gamma \mapsto \nu_{\gamma}$ is continuous (for the topology of the weak convergence).
% and $\nu_{\gamma}$ is absolutely continuous with respect to the Lebesgue measure. 
\end{theorem}

\begin{figure}[htb]
\begin{center}
 \includegraphics[angle=0,width = 10cm, height=5cm]{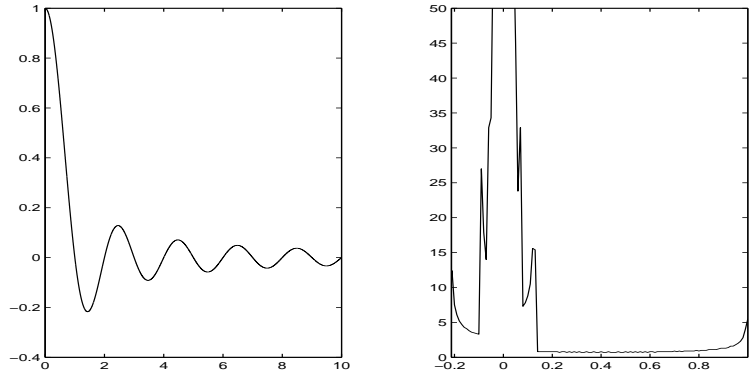}
\caption{Left: $\hat f (\xi)$ for $f(x) = \ind(0 \leq |x| \leq r)$. Right: the level set function $\psi$.}
\label{fig1}
\end{center}\end{figure}

The exact computation of $\nu_\gamma$ is a difficult problem, we
will compute the value $\nu_\gamma (P_m)$, where $P_m$ is the
polynomial $t \mapsto t^m$ (Equation (\ref{eq:momentnu2})). However,
the behavior of $\nu_\gamma$ as $\gamma$ goes to infinity is
simpler. Indeed, we define the Fourier transform of $f$ by, for all
$\xi \in \R^d$, $\hat f(\xi) = \int_0 ^ \infty e^{-2i\pi \xi. x}
f(x)dx$. Since $f$ has a bounded support, $\hat f$ is infinitely differentiable. We assume that the following inversion formula holds
\begin{equation}
\label{hyp:f}
\hbox{a.e. and at $0$,} \quad   f(x) = \int_{\R^d} \hat f(\xi) e^{2i \pi \xi.x} d\xi.
\end{equation}
Note that if $f$ is hermitian ($f(-x) = \overline f(x)$) then $\hat f(\xi) \in \mathbb R$ and for $\epsilon > 0$, $\int_{\R ^d} \ind( |\hat f(\xi)| \geq \epsilon ) d\xi$ is finite. Hence by the change of variable formula,  there exists a function $\psi$ such that  for all continuous functions $h$ with $0 \notin
supp(h)$:
 \begin{equation*}
\label{eq:levelset}
\int_{\R} h(t) \psi(t) dt = \int_{\R^d} h( \hat  f(x)) dx.
\end{equation*}
$\psi$ is the level sets function of $\hat f$, if $\ell$ denotes the $d$-dimensional Lebesgue measure, for all $t > 0$, $\psi (t) = \lim_{\epsilon \to 0} \ell(\{x : |\hat f(x) - t| \leq \epsilon \})/ \epsilon$. If $d = 1$ and $\hat f$ is a diffeomorphism from $\R$ to $K$ then $\psi$ has support on $K$ and is equal to $\psi(t)  = (\hat f ^{-1})'(t) $. 
\begin{theorem}
\label{th:alphamain} If $f$ is hermitian and (\ref{hyp:f}) holds true, then as $\gamma$ goes
to infinity, for all analytic functions $h(t) = \sum_{m \in \N} h_m t^m $ with $h_0 = 0$ and $\sum_{m \in \N} |h_m| t^m$ finite for all $t$:
$$
\int_{\R} h(t) \nu_\gamma (dt) \sim \int_{\R} h(t)  \gamma^{-2} \psi
(\frac t  \gamma)dt = \int_{\R} \gamma^{-1} h(\gamma t)   \psi (t
)dt.
$$
\end{theorem}
This result states that the measure $\nu_\gamma (dt)$ is  in a weak sense equivalent (not in the measure theory sense) to the measure $\gamma^{-2} \psi ( t / \gamma ) dt$ in the high density asymptotic. The measure $\psi (t)dt$ is the continuous analog of the counting measure $\mu$ in Theorem \ref{th:main}. As an example, if $d = 1$ and $f(x) = \ind(0 \leq |x| \leq r)$ then $\hat f(\xi) = r \sinc (2 \pi \xi r)$ and $\psi(t)$ is plotted in Figure \ref{fig1}.

\noindent {\bf Remark.}  
Let $\gamma_n = n \delta_n ^ d$, if $\gamma_n$   tends to infinity and $\delta_n$ goes to $0$, with the material of this note, we may also prove the convergence of $\delta_n ^ d \sum_{i=1}^ n \delta_{\lambda'_i(n)/ \gamma_n}$ to the measure $\psi(t) dt$ on all continuous function $h$ with compact support and $ 0 \notin supp(h)$.

The spectral radius of the matrix $B_n$ is not computed explicitly in this paper. However, the following upper bound is available:
\begin{proposition}
\label{prop:upperbound}
If $ d\geq 2$ and $Po( \gamma)$ denotes a random variable with Poisson distribution of intensity $\gamma$, then with a probability tending to $1$ as $n$ goes to infinity, 
$$
\max_{ 1 \leq i \leq n} |\lambda_i' (n) | \leq  j(n) \sup_{ x \in \Omega} |f(x)|, 
$$ 
where $j(n)$ is solution of:  $n \PP ( Po(\gamma) \geq j(n) +1 ) \leq 1 < n \PP(  Po(\gamma) \geq j(n)).$
\end{proposition}
For $n$ large enough, using the inequality $\PP ( Po (\gamma) \geq k ) \leq \exp ( - \frac k 2 \ln ( \frac k \gamma) )$, for $k \geq e^2 \gamma$, we deduce that $j(n) \leq 3 \ln n / \ln \ln n$ for $n$ large enough.

The remainder of this paper is organized as follows. In Section \ref{sec:1}, we prove Theorem \ref{th:main}, In Section \ref{sec:2}, we prove Theorems \ref{th:scmain}, \ref{th:alphamain} and Proposition \ref{prop:upperbound}. Finally, in Section \ref{sec:misc}, we state some simple results on the eigenvectors of $A$ and on the correlation of the eigenvalues.

By convention $C$ will denote a constant which does not depend on $n$. Its exact value may change throughout the paper. Also we define: $\| F\|_{\infty} = \sup_{ x \in  \R^d} |F(x)|$ and $B(x,r)$ will denote the open ball of radius $r$ and center $x$ on the torus $\T^d$.

\section{Proof of Theorem \ref{th:main}}

\label{sec:1}

The proof of Theorem \ref{th:main} relies on the classical Wigner's method \cite{wigner} to compute the empirical mean distribution measure of eigenvalues. We will compute for all $m \in \N$:
$$
\EE \tr A_n ^m = \frac {1} {n ^m} \EE \tr A ^m = \frac {1} {n ^m}  \sum_{ i=1} ^n \lambda_i ^m = \mu_n (P_m).
$$
We will then use a Talagrand's concentration inequality to prove that $\tr A_n ^m$ is not far from its mean and conclude. About the rate of convergence of $\mu_n$ to $\mu$, we will
state (in the forthcoming Lemma \ref{le:wigner}) that, if $P_m (t) =
t^m$, $m\geq 1$,
\begin{equation}
\label{eq:weakmu2} \lim_n n \Bigm( \EE \mu_n ( P_m ) - \mu(P_m)\Bigm) =
\sum_{q = 1} ^ {m-1}q  \mu(P_q) \mu(P_{m-q}) -  \frac{m(m-1)}{2}
\mu(P_m).
\end{equation}

We begin with a technical lemma.
\begin{lemma}
\label{le:combi} For $0 \leq p \leq m$, let $\Sigma_{m,p}$ be the
set of surjective mappings from $\{1,\cdots,m\}$ to
$\{1,\cdots,p\}$. We have:
\begin{equation}
\label{eq:combi} \EE \tr A^m = \sum_{p = 1}^{m} \binom{n}{p}
\sum_{\phi \in \Sigma_{m,p}} \int_{\Omega^{p}} \prod_{j=1}^ m
F(x_{\phi(j)} - x_{\phi(j+1)}) dx_1 \cdots dx_{p},
\end{equation}
with $\phi(m+1) = \phi(1)$ and with the convention that
$\binom{n}{p} = 0$ for $p > n$.
\end{lemma}
\bp
By definition:
\begin{equation}
\label{eq:trAm}
\tr A^ m = \sum_{i_1 ,\cdots i_m} \prod_{j=1} ^ m F(X_{i_j} - X_{i_{j+1}}),
\end{equation}
with $i_{m+1} = i_1$ and the sum is over all n-tuples of integers
${\bf i}  = (i_1,\cdots,i_m)$ in $\{1,n\}^m$. Let $p (\bf i)$ be the set of distinct indices in $\bf i$. We can
define a surjective mapping $\phi_{\bf i}$ in $\Sigma_{m,p(\bf i)}$
such that $i_j = i_{\phi_{\bf i}(j)}$. Taking the expectation in
Equation (\ref{eq:trAm}), we get
$$
\EE \tr A^m = \sum_{{\bf i} = (i_1 ,\cdots i_m) } \int_{\Omega^m}
\prod_{j=1}^ m F(x_{\phi_{\bf i}(j)} - x_{\phi_{\bf i}(j+1)}) dx_1
\cdots dx_{p(\bf i)},
$$
We then reorder the terms. We consider the equivalence relation in
$\Sigma_{m,p}$, $\phi \sim \phi'$ if there exists a permutation
$\sigma$ of $\{1, \cdots,p\}$ such that $\sigma\circ \phi =
\phi'$. The value of $ \int_{\Omega^m} \prod_{j=1}^ m
F(x_{\phi(j)} - x_{\phi(j+1)})dx_1 \cdots dx_{p}$ is constant on
each equivalence class. Let $\phi \in \Sigma_{m,p}$, the numbers of
indices $\bf i$ such that $\phi_{\bf i} \sim \phi $ is equal to $n!
/ (n-p)!$ (if $n \geq p$ and $0$ otherwise). Since there are $p!$
surjective mappings in the class of equivalence of $\phi$, we deduce
Equation (\ref{eq:combi}).\ep

\begin{lemma}
\label{le:wigner}
For each $m$,
$$
\EE \mu_n (P_m) = \mu(P_m) + \frac{1} {n}  \Bigm( \sum_{q = 1}
^ {m-1} q \mu(P_q) \mu(P_{m-q}) -  \frac{m(m-1)}{2} \mu(P_m)\Bigm) +
o(\frac 1 n).
$$
\end{lemma}

\bp We apply Lemma \ref{le:combi} and identify the  coefficients in
$n^m$ and $n^{m-1}$ in Equation (\ref{eq:combi}). We first consider
the term in $n^m$, such a term comes from $p=m$:
$$
\frac{n !}{(n-m)!}  \int_{\Omega^m} \prod_{j=1}^ m F(x_{j} - x_{j+1}) dx_1 \cdots dx_{m},
$$
By induction, we easily obtain that
$$
\int_{\Omega^m} \prod_{j=1}^ m F(x_{j} - x_{j+1}) dx_1 \cdots dx_{m} = \int_\Omega F^{*m} (0) dx_1 = F^{*m} (0),
$$
where $*$ denotes the convolution operator: $F * G (y) = \int_{\Omega} F(y-x) G(x) dy$ and $F^{*m}$ is $F* F \cdots *F$ ($m$ times). We recall the two properties: $\widehat {F * G} (k) = \hat F (k) \hat G(k)$ and  $
F(0)  = \sum_{ k \in \Z^d} \hat F(k),
$
in order to get:
$$
\int_{\Omega^m} \prod_{j=1}^ m F(x_{j} - x_{j+1}) dx_1 \cdots dx_{m} = \sum_{k \in \Z^d} \hat F(k)^ m = \mu(P_m).
$$
We thus deduce that:
$$
\lim_n \EE \mu_n (P_m) = \mu(P_m).
$$
It remains to identify the terms in $n^{m-1}$ in Equation
(\ref{eq:combi}). This term comes from two contributions $p = m$ and
$p=m-1$. Since $n ! / (n-m)!  = n^m - n^{m-1} \sum_{i=0}^{m-1} i +
o(n^{m-1})$, the term in $n^{m-1}$ in $p=m$ is equal to:
\begin{equation}
\label{eq:wigner1}
- \frac{m(m-1)}{2} \mu(P_m).
\end{equation}
The leading term for $p=m-1$ is
\begin{equation}
\label{eq:wigner2}
\frac{n!}{(n-m+1)!(m-1)!}  \sum_{\phi \in \Sigma_{m,m-1}} \int_{\Omega^{m-1}} \prod_{j=1}^ m F(x_{\phi(j)} - x_{\phi(j+1)}) dx_1 \cdots dx_{m-1},
\end{equation}
Now if $\phi \in \Sigma_{m,m-1}$, $\phi^{-1} (i)$ is not reduced  to a
single point for a unique index $i_\phi$. Since the value of $
\int_{\Omega^{m-1}}\prod_{j=1}^ m F(x_{\phi(j)} - x_{\phi(j+1)})
dx_1 \cdots dx_{m-1}$ is invariant under permutations of the
indices, without loss of generality, we may assume that $i_\phi =
1$, $\phi^{-1} (1) = \{1,q+1\}$ with $q \in \{1,\cdots, m-1\}$ and $\phi(j)
= j$ if $j \leq q$ and $\phi(j) = j+1$ if $j > q+1$. For such
$\phi$,  integrating over $x_2, \cdots, x_q, x_{q+1}, \cdots ,x_{m-1}$,
\begin{eqnarray*}
\int_{\Omega^m} \hspace{-3pt} \prod_{j=1}^ m F(x_{\phi(j)} - x_{\phi(j+1)}) dx_1  \cdots  dx_{m-1} 
%& = &  \int_{\Omega^m} \hspace{-8pt} F(x_{1} - x_{2}) \cdots  F(x_q - x_1) \\
%&&  \quad \quad \quad   F(x_1 - x_{q+1}) \cdots F(x_{m-1} - x_1) dx_1 \cdots dx_{m-1} \\
&  = & \int_\Omega F^{*(q)} (0) F^{*(m-q)} (0) dx_1 \\
& = &  \mu(P_q) \mu(P_{m-q}).
\end{eqnarray*}
Finally, for each $q$, there are $(m-q) \times (m-1) ! $ surjective
mappings such that, up to a permutation of the indices, $\phi^{-1}
(1) = \{1,q+1\}$ and $\phi(j) = j$ if $j \leq q$ and $\phi(j) = j+1$
if $j > q+1$.  Indeed for such $\phi$, there are $(m-q)$ possible pairs $(i_1,i_1+q)$, $1 \leq i_1 
\leq m - q $ such that $\phi(i_1) = \phi(i_1 + q)$.  Therefore Equation (\ref{eq:wigner2}) can be
written as:
\begin{eqnarray}
\label{eq:lemmap2} \frac{n! (m-1)! }{(n-m+1)!(m-1)!}  \sum_{q =
1} ^ {m-1} (m-q) \mu(P_q) \mu(P_{m-q}) = n^{m-1} \sum_{q = 1} ^
{m-1} q  \mu(P_q) \mu(P_{m-q}) + o (n^{m-1}).
\end{eqnarray}
Adding this last term with the term (\ref{eq:wigner1}), we get the stated formula.\ep

We may now prove Theorem \ref{th:main}.

\noindent{\it Proof of Theorem \ref{th:main}.}

We fix $n$ and for each $m \geq 1$, we define the functional:
$$
Q_m({\cX_n}) = \frac 1 {n^{m-1}} \tr A^m = n \mu_n (P_m).
$$
If ${\bf x} ,{\bf y} \in \Omega^n$, let $d({\bf x} ,{\bf y}) = \sum_{i=1} ^ n \ind( x_i \neq y_i)$ denote the Hamming distance.
The functional $Q_m$ is Lipschitz for the Hamming distance $d$. Indeed, define ${\bf x^l} = ( x^l_j)_{1 \leq j \leq n}$ by $x^l_j = x_j$ for $j \ne l$ and $ x^l_l \ne x_l$, we have:
\begin{eqnarray*}
\Bigm|Q_m({\bf x}) - Q_m({\bf x^l})\Bigm| & = & \frac{1} { n^{m-1}} \Bigm| \sum_{i_1,\cdots,i_m}  \prod_{j=1} ^ m F(x_{i_j} - x_{i_{j+1}}) - \prod_{j=1} ^ m F(x^l_{i_j} - x^l_{i_{j+1}}) \Bigm| \\
& \leq & 2 m \| f\|^m_\infty,
\end{eqnarray*}
indeed $|\prod_{j=1} ^ m F(x_{i_j} - x_{i_{j+1}}) - \prod_{j=1} ^ m F(x^l_{i_j} - x^l_{i_{j+1}})|$ is at most $2 \| f\|^m_\infty$ and it is non zero only if there exists a index $i_j$ such that $i_j = l$. It follows easily that $Q_m$ is $2 m \| f\|^m_\infty$-Lipschitz for the Hamming distance $d$.

Let $M_m$ denote the median of $Q_m$. We may apply a Talagrand's Concentration Inequality (see for example Proposition 2.1 of Talagrand \cite{talagrand96}),
$$
\PP ( | Q_m - M_m | > t) \leq 4 \exp(- \frac{ t^2 } {4 m^2 \| f\|^{2m}_\infty n}),
$$
integrating over all $t$ we deduce:
$$
| n \EE \mu(P_m) - M_m | \leq \EE | \mu(P_m) - M_m | \leq C_m \sqrt n,
$$
for some constant $C_m$ and it follows, that for all $s > C_m  / \sqrt n$:
$$
\PP ( | \mu_n (P_m)  - \EE \mu(P_m) | > s) \leq 4 \exp(- n \frac{ (s -C_m  / \sqrt n ) ^2 } {4 m^2 \| f\|^{2m}_\infty}),
$$
Using the Borel Cantelli Lemma and Lemma \ref{le:wigner}, a.s.
$
\lim_n \mu_n (P_m) = \mu(P_m).
$
\ep

\section{Limit Spectral Measure of Scaled ERM}
\label{sec:2}

\subsection{Proof of Theorem \ref{th:scmain}}

The study of the first model was simplified by the absence of boundary effects with $\Omega$. So in order to prove Theorem \ref{th:scmain}, we will first discard them in the second model. We define $F_\delta$ as the 1-periodic extension of $f_{\delta}$: for all $x \in \R^d$, there exists a unique couple $(y,u)$ such that $x = y + u$, with $u \in \Z^d$ and $y \in  \Omega$, and we set $F_{\delta} (x) = f_\delta (y)$.

We now introduce a matrix and its spectral empirical measure: 
$$\tilde B_n = (F_{\delta_n}
(X_i - X_j))_{ 1 \leq i \leq j \leq n} \quad \hbox{ and }
 \quad \tilde \nu_n = \frac 1 n \sum_{i=1} ^n \delta_{\tilde \lambda_i(n)},
$$
where  $(\tilde \lambda_1(n), \cdots, \tilde \lambda_n (n) )$ is the spectrum of $\tilde B_n$. The next lemma states that the limiting spectral measures of $\tilde \nu_n$ and $\nu_n$ are equal.
\begin{lemma}
\label{le:boundary}
For the topology of the weak convergence of (signed) measures, a.s.
$ \nu_n - \tilde \nu_n $ converges as $n$ goes to infinity to the null measure.
\end{lemma}

\bp
It is sufficient to prove that for all $m \geq 1$, a.s.  $\lim_n \nu_n (P_m) - \tilde \nu_n (P_m) = 0$.  To this end, we notice that if $x,y \in \Omega$, $f_{\delta} ( x - y ) = F_{\delta} (x-y)$ unless $x \in \Omega \backslash (1- \delta) \Omega$ and $y \in B(x, \delta)$.  We write:
\begin{eqnarray*}
\bigm| \nu_n (P_m) - \tilde \nu_n (P_m) \bigm|
 & \leq & \frac 1 n \sum_{i_1, \cdots,i_m} \bigm|  \prod_{j=1} ^ m f^m_{\delta_n} (X_{i_{j}} - X_{i_{j+1}}) -   \prod_{j=1} ^ m F_{\delta_n} (X_{i_{j}} - X_{i_{j+1}}) \bigm|   \\
& \leq & \frac 1 n \sum_{i_1, \cdots,i_m} 2 \|f\|^m_{\infty} \ind( X_{i_1} \in \Omega \backslash (1- \delta_{n}) \Omega) \prod_{j=2}  ^ m \ind( X_{i_j} \in B(X_{i_1}, m \delta_n)) \\
%& \leq & \frac {2} {n}  \|f\|^m_{\infty}  \sum_{i=1}^ n  \ind( X_i \in  \Omega \backslash (1- \delta_{n}) \Omega) \frac{N_n(B(X_i,m\delta_n))^m}{m^m},
& \leq & \frac {2} {n}  \|f\|^m_{\infty}  N_n(\Omega \backslash (1- m \delta_n) \Omega)^m,
\end{eqnarray*}
where $N_n$ is the counting measure $N_n (\cdot) = \# \{i \in \{1, \cdots,n\}  : X_i \in \cdot \}$.
%Note that $\EE N_n (B(X_i,m\delta_n) )^m \leq C_m (n \delta_n ^ d)^m \leq C$ and $\EE N_n( \Omega \backslash (1- \delta_{n}) \Omega) / n \leq C \delta_n$. It follows  classically that $ 1/n \sum_{i=1}^ n  \ind( X_i \in  \Omega \backslash (1- \delta_{n}) \Omega) N_n(B(X_i,m\delta_n))^m$ tends a.s. to $0$. 
Note that $\PP ( X_1 \in \Omega \backslash (1- m \delta_n) \Omega ) \leq C \delta_n$. By the strong law of large numbers, it follows easily that that $N_n(\Omega \backslash (1- m \delta_n) \Omega) / n$ converges almost surely to $0$. \ep

By Lemma \ref{le:boundary}, we may focus on $\tilde B_n$ and $\tilde \nu_n$. In order to keep the notations as light as possible we drop the '' $\tilde \cdot$ '' in $\tilde B_n$ and $\tilde \nu_n$.

We first prove that,
\begin{equation}
\label{eq:step1}
 \nu_n   \hbox{ converges in probability to a measure } \nu_{\gamma} \hbox{ for the weak convergence.}
\end{equation}
By Lemma \ref{le:combi}, if $m \geq 1$,
\begin{equation}
\label{eq:scwig} \EE \nu_n (P_m) =   \frac 1 n \sum_{p = 1}^{m}
\binom{n}{p} \sum_{\phi \in \Sigma_{m,p}} \int_{\Omega^{p}}
\prod_{j=1}^ m F_{\delta_n} (x_{\phi(j)} - x_{\phi(j+1)}) dx_1
\cdots dx_{p}. \end{equation} We begin with an elementary lemma.
\begin{lemma}
\label{le:varaveugle} If $\phi \in \Sigma_{m,p}$, $p>1$ the value of
$$
\int_{\Omega^{p-1}}  \prod_{j=1}^ m F (x_{\phi(j)} - x_{\phi(j+1)})
dx_2 \cdots dx_{p}
$$
does not depend on $x_{1}$.
\end{lemma}

\bp We consider the change of variable, for $j > 1$, $x'_j = x_j -
x_{1}$. The Jacobian of this change of variable is $1$. If we set
$x'_{1} = 0$, we obtain $\int_{\Omega^{p-1}}  \prod_{j=1}^ m F
(x_{\phi(j)} - x_{\phi(j+1)}) dx_2 \cdots dx_{p} =
\int_{\Omega^{p-1}}  \prod_{j=1}^ m F (x'_{\phi(j)} -
x'_{\phi(j+1)}) dx'_2 \cdots dx'_{p}$. \ep

Assume $m \geq 2$, by Lemma \ref{le:varaveugle}, we have:
\begin{eqnarray}
\EE \nu_n (P_m)  & = &   F_{\delta_n}(0) ^ m  + \frac 1 n \sum_{p = 2}^{m} \binom{n}{p} \sum_{\phi \in \Sigma_{m,p}} \int_{\Omega^{p-1}} \prod_{j=1}^ m F_{\delta_n} (x_{\phi(j)} - x_{\phi(j+1)}) dx_2 \cdots dx_{p} \nonumber \\
& = & f(0)^ m + \frac 1 n \sum_{p = 2}^{m} \binom{n}{p}\hspace{-3pt} \sum_{\phi
\in \Sigma_{m,p}} \hspace{-3pt} \Delta (\phi) + \int_{\Omega^{p-1}} \hspace{-3pt}\prod_{j=1}^
m f_{\delta_n} (x_{\phi(j)} - x_{\phi(j+1)}) dx_2 \cdots dx_{p}
\label{eq:temp1sc}
\end{eqnarray}
where $\Delta (\phi)  = \int_{\Omega^{p-1}} \prod_{j=1}^ m
F_{\delta_n} (x_{\phi(j)} - x_{\phi(j+1)})- \prod_{j=1}^ m f_{\delta_n} (x_{\phi(j)} -
x_{\phi(j+1)}) dx_2 \cdots dx_{p}$. Since the support of
$f_{\delta}$ is included in $\delta \Omega$, if $f_{\delta}
(x_{\phi(j)} - x_{\phi(j+1)}) \ne F_{\delta} (x_{\phi(j)} -
x_{\phi(j+1)}) $ then $x_{\phi(j)},  x_{\phi(j+1)} \in \Omega
\backslash (1-\delta) \Omega$. Moreover notice that if $
\prod_{j=1}^ m F (x_{\phi(j)} - x_{\phi(j+1)}) \neq 0$ then
$x_2,\cdots x_{p} \in B(x_{1},(m-1)\delta)$.  By Lemma \ref{le:varaveugle}, from
now on, we can assume without loss of generality:$$x_{1} = 0,$$ and
then $\Delta(\phi) = 0$ for $\delta < 1/(2m)$.

Considering the change of variable $y_i = x_i / \delta_n$ in the integrands of Equation (\ref{eq:temp1sc}), we obtain, for $\delta < 1/(2m)$, with $y_1 = 0$,
\begin{eqnarray}\label{eq:deltasc2}
\EE \nu_n (P_m) =  f(0)^ m   +  \sum_{p = 2}^{m} \frac
{\delta_n^{d(p-1)}}{ n} \binom{n}{p} \sum_{\phi \in \Sigma_{m,p}}
\int_{(\delta_n ^{-1} \Omega)^{p-1}} \prod_{j=1}^ m f (y_{\phi(j)} -
y_{\phi(j+1)}) dy_2 \cdots dy_{p}.
\end{eqnarray}

Finally, since $\binom{n}{p} \sim n^{p} / p!$ as $n$ goes to
infinity, we deduce that, for $m \geq 2$,
\begin{equation}
\label{eq:momentnu} \lim_{n \to \infty}  \EE \nu_n (P_m) =   f(0)^ m
+  \sum_{p = 2}^{m} \frac{ \gamma^{p-1}}{p!} \sum_{\phi \in
\Sigma_{m,p}} \int_{(\R^d)^{p-1}} \prod_{j=1}^ m f (y_{\phi(j)} -
y_{\phi(j+1)}) dy_2 \cdots dy_{p}.
\end{equation}
(For $m \leq 1$, we have $\nu_n (P_0) = 1$ and $\nu_n(P_1) = f(0)$).

We check easily that the right hand side of Equation (\ref{eq:momentnu}) is bounded by $(C m)^{m}$ for some constant $C$ not depending on $m$. Therefore, by Carleman's Condition, there is exists a unique measure $\nu_{\gamma}$ such that $\lim_{n \to \infty}  \EE \nu_n (P_m) = \nu_{\gamma} (P_m) $. In particular, the sequence $(\nu_n)_{n \in \mathbb N}$ is tight and we have proved (\ref{eq:step1}).

The continuity of $\gamma \mapsto \nu_\gamma$ follows from the comtinuity of $\gamma \mapsto \nu_{\gamma} (P_m)$. Indeed, let $(\gamma_n)_{ n \in \N}$ be a sequence converging to $\gamma < \infty$. Since $\sup_{n} \nu_{\gamma_n} (P_2) < \infty$, the sequence $(\nu_{\gamma_n})_{ n \in \N}$ is tight. Hence for all $\epsilon > 0$, there exists  a compact set $K$ such that for all $n$ $\nu_{\gamma_n} (K^c) \leq \epsilon$.  Now, let $h$ be a continuous function with compact support, we need to prove that $\lim_{n \to \infty} \nu_{ \gamma_n} (h) = \nu_{\gamma} (h)$. Fix $\epsilon$, there exists a polynomial $P$ such that $\sup_{ x \in K} | h(x) - P(x) | \leq \epsilon$, we deduce that $|\nu{\gamma_n} ( h ) - \nu_{\gamma} ( h) | \leq |\nu_{\gamma_n} (h)  - \nu_{\gamma_n} (P) |  +  |\nu_{\gamma_n} (P)  - \nu_{\gamma} (P) |  +  |\nu_{\gamma} (P)  - \nu_{\gamma} (h) | \leq   2\epsilon (1 + \| h\|_{\infty}) +  |\nu_{\gamma_n} (P)  - \nu_{\gamma} (P) |$. Letting $n$ tends to infinity, since $\epsilon$ is arbitrary small and $\gamma \mapsto \nu_{\gamma} (P)$ is continuous, we obtain: $\lim_{n \to \infty} \nu_{ \gamma_n} (h) = \nu_{\gamma} (h)$.

It remains to prove the almost sure convergence of $\nu_n$. We will prove that for each $m \geq 1$, there exists a constant $C$ and
\begin{equation}
\label{eq:4moment}
\EE \bigm(\tr B_n ^ m - \EE \tr B_n ^ m \bigm)^4 \leq C n^2.
\end{equation}
This last equation implies $\EE \bigm( \nu_n(P_m) - \EE \nu_n(P_m)\bigm)^4 \leq C / n^2$ and by Borel Cantelli Lemma, we deduce that $\nu_n(P_m)$ converges almost surely toward $\nu_{\gamma} (P_m)$.

It remains to prove Inequality (\ref{eq:4moment}). A circuit in  $\{1,\cdots,n\}$ of length $m$ is a mapping $\pi : \Z  \to \{1,\cdots n\}$ such that for all integer $r$, $\pi(m+r)  = \pi(r)$. Following Bryc, Dembo and Jiang \cite{dembo2006}, we introduce the new notation:
$$F_{\pi} = \prod_{i=1}^ m F_{\delta_n} ( X_{\pi(i)} - X_{\pi(i+1)}).$$
We then write:
\begin{eqnarray}
\label{eq:trace4}
\EE \bigm(\tr B_n ^ m - \EE \tr B_n ^ m \bigm)^4  =  \EE \bigm( \sum_{\pi}  F_{\pi} - \EE F_{\pi} \bigm)^4
 =  \sum_{\pi_1,\cdots,\pi_4} \EE \bigm[ \prod_{l=1}^4 F_{\pi_l} - \EE F_{\pi_l}\bigm],
\end{eqnarray}
where the sums are over all circuits in  $\{1,\cdots,n\}$ of length $m$.

Notice that $\EE \bigm[ \prod_{l=1}^4 F_{\pi_l} - \EE F_{\pi_l}\bigm] = 0$ if there exists a circuit $\pi_k$, $1 \leq k \leq 4$ such that the image of $\pi_k$ has an empty intersection with the union of the images of $\pi_l$, $l \ne k$. Indeed, due to the independence of the variables $(X_i)_{1 \leq i \leq n}$, $F_{\pi_k} - \EE F_{\pi_k}$ is then independent of  $\prod_{l\ne k} F_{\pi_l} - \EE F_{\pi_l}$.

Two circuits $\pi_1$ and $\pi_2$ in  $\{1,\cdots,n\}$ of length $m_1$ and $m_2$ with a non empty intersection of their images may be concatenated into a circuit in  $\{1,\cdots,n\}$ of length $m_1+m_2$ as follows. Assume that $\pi_1 (i_0) = \pi_2 (j_0)$, we define the circuit $\pi_{1,2}$  of length $m_1+m_2$ by for $i \in \{1, \cdots, m_1 + m_2\}$
$$
\pi_{1}.\pi_2 (i) = \left\{ \begin{array}{lcl} \pi_1(i) & \hbox{if} & 1 \leq i \leq i_0 \\
\pi_2(i-i_0+j_0) & \hbox{if} & i_0 + 1 \leq i \leq i_0 +m_2 \\
\pi_1(i-m-i_0) & \hbox{if} & i_0 + m +1 \leq i \leq m_1 + m_2
\end{array}\right.
$$
We have:
$$
F_{\pi_1} F_{\pi_2} = F_{\pi_1 \cdot \pi_2}.
$$
Using the same reasoning as for Equation (\ref{eq:deltasc2}), we get
  $$\EE F_{\pi_1} F_{\pi_2} =  \delta_n^{d(q-1)} \int_{(\delta_n ^{-1} \Omega)^{q-1}} \prod_{j=1}^ {2m} f (y_{\pi_1.\pi_2(j)} - y_{\pi_1.\pi_2(j+1)}) dy_{i_2} \cdots dy_{i_q},$$
where is $q = q(\pi_1,\pi_2)$ is the cardinal of the union of the images of $\pi_{1}$ and $\pi_2$ and $(y_{i_1},\cdots,y_{i_q})$ is the image of  $\pi_{1} \cdot\pi_2$ and $y_{i_1} = 0$.

If $N(\pi_1,\pi_2) $ is the cardinal of the intersection of the images of $\pi_1$ and $\pi_2$, if $N(\pi_1,\pi_2)\geq 1$, we obtain
\begin{equation}\label{eq:2termes}
\EE | F_{\pi_1} F_{\pi_2} | \leq C  n ^{- q(\pi_1,\pi_2) +1}.
\end{equation}
Otherwise, $N(\pi_1,\pi_2) = 0 $, if $q(\pi_i)$ is the cardinal of of the image of $\pi_i$,
\begin{equation}
\EE | F_{\pi_1} F_{\pi_2}| = \EE| F_{\pi_1}| \EE | F_{\pi_2}| \leq  C  n ^{- q(\pi_1) - q(\pi_2)+2} =   C  n ^{- q(\pi_1,\pi_2)+2}
\label{eq:prod2}.
\end{equation}

Similarly assume that $N(\pi_1,\pi_2) \geq 1$, $N(\pi_1.\pi_2 , \pi_3) \geq 1$,  if $q(\pi_1,\pi_2,\pi_3)$ is the cardinal of the union if the images of $\pi_1, \pi_2 , \pi_3$ then $F_{\pi_1} F_{\pi_2} F_{\pi_3} = F_{(\pi_1.\pi_2 ).\pi_3}$ and we deduce similarly $$
\EE| F_{\pi_1} F_{\pi_2}  F_{\pi_3} | \leq C n ^{- q(\pi_1,\pi_2,\pi_3) +1}.
$$
Finally assume that $N(\pi_1,\pi_2) \geq 1$, $N(\pi_1.\pi_2 , \pi_3) \geq 1$,  $N(\pi_1.\pi_2.\pi_3 , \pi_4) \geq 1$, if  $q(\pi_1,\pi_2,\pi_3)$ is the cardinal of the union if the images of $\pi_1, \pi_2 , \pi_3,\pi_4$, we obtain:
\begin{equation}
\label{eq:4termes}
\EE |F_{\pi_1} F_{\pi_2}  F_{\pi_3} F_{\pi_4} | \leq C n ^{- q(\pi_1,\pi_2,\pi_3,\pi_4) +1}.
\end{equation}

By (\ref{eq:trace4}), it remains to decompose:
$$
4 ! \sum_{(\pi_1,\cdots,\pi_4)\in S \cup S'} \EE \bigm[ \prod_{l=1}^4 F_{\pi_l} - \EE F_{\pi_l}\bigm]
$$
where $S$ is the set of quadruples of circuits such that $N(\pi_1,\pi_2) \geq 1$, $N(\pi_1.\pi_2 , \pi_3) \geq 1$,  $N(\pi_1.\pi_2.\pi_3 , \pi_4) \geq 1$ and $S'$ is the set of quadruples of circuits such that $N(\pi_1,\pi_2) \geq 1$ and $N(\pi_3,\pi_4) \geq 1$ and otherwise for $i < j$, $N(\pi_i,\pi_j) = 0$.

The decomposition of the $\EE \bigm[ \prod_{l=1}^4 F_{\pi_l} - \EE F_{\pi_l}\bigm]$ gives rise to four types of terms:
\begin{enumerate}
\item $\sum_{(\pi_1,\cdots,\pi_4) \in S \cup S'}  \prod_{l=1} ^ 4 \EE F_{\pi_l}$,
\item $\sum_{(\pi_1,\cdots,\pi_4)\in S \cup S'}   \EE \prod_{l=1} ^ 4  F_{\pi_l}$,
\item $\sum_{(\pi_1,\cdots,\pi_4) \in S \cup S'} \EE F_{\pi_{l_1}}  F_{\pi_{l_2}} \EE F_{\pi_{l_3}} F_{\pi_{l_4}}$,
\item $\sum_{(\pi_1,\cdots,\pi_4) \in S \cup S'} \EE F_{\pi_{l_1}}  \EE  F_{\pi_{l_2}}  F_{\pi_{l_3}} F_{\pi_{l_4}}$,
\end{enumerate}
where $(l_1, l_2, l_3,l_4)$ is a permutation of $(1,2,3,4)$. We will apply successively the same method to bound these terms.

We begin with the terms of type 1, we have: $\prod_{l=1} ^ 4 \EE F_{\pi_l} \leq C n^{-\sum_{l=1} ^ 4 q(\pi_l)+4}$.
Since $(\pi_1,\cdots,\pi_4) \in S \cup S'$, $q(\pi_1,\pi_2,\pi_3, \pi_4)  \leq \sum_{l=1} ^ 4 q(\pi_l)-2$, hence:
$$
\prod_{l=1} ^ 4 \EE F_{\pi_l} \leq C n^{- q(\pi_1,\pi_2,\pi_3,\pi_4)+2}.
$$
There are at most $C n^{q}$ quadruples of circuits such that $q(\pi_1,\pi_2,\pi_3, \pi_4) = q$, therefore the terms of type $1$ may be bounded as by
 $$\sum_{(\pi_1,\cdots,\pi_4) \in S \cup S'}  \prod_{l=1} ^ 4 \EE | F_{\pi_l}| \leq C n ^2.$$

We now deal with the terms of type $2$. By (\ref{eq:4termes}), if $(\pi_1,\cdots,\pi_4) \in S$, $ \EE \prod_{l=1} ^ 4  | F_{\pi_l} |\leq C n ^{- q(\pi_1,\pi_2,\pi_3,\pi_4) +1}$ otherwise $(\pi_1,\cdots,\pi_4) \in S'$ and, by (\ref{eq:prod2}),  $ \EE \prod_{l=1} ^ 4  | F_{\pi_l} |\leq C n ^{- q(\pi_1,\pi_2,\pi_3,\pi_4) +2}$. There are at most $C n^{q}$ quadruples of mappings such that $q(\pi_1,\pi_2,\pi_3, \pi_4) = q$. Hence
$$\sum_{(\pi_1,\cdots,\pi_4)\in S \cup S'} \EE  \prod_{l=1} ^ 4  | F_{\pi_l}| \leq C n^2 .$$

We turn to the terms of type $3$: $\sum_{(\pi_1,\cdots,\pi_4) S \cup S'} \EE F_{\pi_{l_1}}  F_{\pi_{l_2}} \EE F_{\pi_{l_3}} F_{\pi_{l_4}}$. Assume first that the quadruple $(\pi_1,\cdots,\pi_4) \in S'$. If $l_1 = 1$, $l_2 = 3$, $l_3 = 2$, $l_4 = 4$, then $\EE F_{\pi_{l_1}}  F_{\pi_{l_2}} \EE F_{\pi_{l_3}} F_{\pi_{l_4}}  = \prod_{l=1} ^ 4 \EE F_{\pi_l}$ and we obtain the same bound that the terms of type $1$. The other cases reduce to the case $l_1 = 1$, $l_2 = 2$, $l_3 = 3$, $l_4 = 4$ and by (\ref{eq:2termes}),  $\EE F_{\pi_{1}}  F_{\pi_{2}} \leq Cn^{-q(\pi_1,\pi_2) +1}$.  There are at most $Cn^{q+q'}$ quadruples such that  $q(\pi_1,\pi_2) = q$ and $q(\pi_3,\pi_4) = q'$.  We deduce  that $\sum_{(\pi_1,\cdots,\pi_4) \in S'} \EE F_{\pi_{1}}  F_{\pi_{2}} \EE F_{\pi_{3}} F_{\pi_{4}} \leq C n^2$.

Assume now that that  $(\pi_1,\cdots,\pi_4) \in S$. We have: $$\EE F_{\pi_{l_1}}  F_{\pi_{l_2}} \EE F_{\pi_{l_3}} F_{\pi_{l_4}} \leq C n^{-q(\pi_{l_1},\pi_{l_2}) - q(\pi_{l_3}, \pi_{l_4}) +2 + \ind(N(\pi_{l_1},\pi_{l_2}) =  0) + \ind(N(\pi_{l_3},\pi_{l_4}) =  0)}.$$ If $N(\pi_{l_1},\pi_{l_2}) =  0$, then $N(\pi_{l_3},\pi_{l_4}) \geq 1 $  and there are at most $Cn^{q+q'-2}$ quadruples such that  $q(\pi_{l_1},\pi_{l_2}) = q$ and $q(\pi_{l_3},\pi_{l_4}) = q'$. Indeed, since $(\pi_1,\cdots,\pi_4) \in S$, the cardinal of the intersection of the images of $ (\pi_{l_1},\pi_{l_2})$  and $(\pi_{l_3},\pi_{l_4})$ is at least $2$. The other cases reduce to the case, $N(\pi_{l_1},\pi_{l_2}) \geq 1$ and $N(\pi_{l_3},\pi_{l_4}) \geq 1 $, for such cases, we notice that there are at most $Cn^{q+q'}$ quadruples such that  $q(\pi_{l_1},\pi_{l_2}) = q$ and $q(\pi_{l_3},\pi_{l_4}) = q'$. In all cases, we conclude that: $\sum_{(\pi_1,\cdots,\pi_4)\in  S} \EE F_{\pi_{1}}  F_{\pi_{2}} \EE F_{\pi_{3}} F_{\pi_{4}} \leq C n^2$.
Hence,
$$\sum_{(\pi_1,\cdots,\pi_4) \in S \cup S'} \EE F_{\pi_{l_1}}  F_{\pi_{l_2}} \EE F_{\pi_{l_3}} F_{\pi_{l_4}} \leq Cn^{2}.$$

It remains to treat the terms of type $4$.
%Assume first that $(\pi_1,\cdots,\pi_4) \in S'$ then all cases reduce to the case: $l_i = i$, $ 1 \leq i \leq 4$. We obtain
%$$\EE F_{\pi_{1}} \EE  F_{\pi_{2}} F_{\pi_{3}} F_{\pi_{4}} \leq C n^{-q(\pi_1) - q(\pi_2) - q(\pi_3,\pi_4) +3}.$$
%Since $(\pi_1,\cdots,\pi_4) \in S'$, $q(\pi_1) + q(\pi_2) \geq q(\pi_1,\pi_2) +1$ and $q(\pi_1,\pi_2) + q(\pi_3,\pi_4) = q(\pi_1,\pi_2,\pi_3,\pi_4)$. We thus have $\EE F_{\pi_{1}} \EE  F_{\pi_{2}} \EE F_{\pi_{3}} F_{\pi_{4}} \leq C n^{-  q(\pi_1,\pi_2,\pi_3,\pi_4)+2}$ and it follows
%$$\sum_{(\pi_1,\cdots,\pi_4)\in   S'} \EE F_{\pi_{1}}  \EE  F_{\pi_{2}}  F_{\pi_{3}} F_{\pi_{4}} \leq Cn^2.$$
Assume that  $(\pi_1,\cdots,\pi_4) \in S \cup S'$, we have:
 \begin{equation}
\label{eq:type3}
\EE F_{\pi_{l_1}} \EE  F_{\pi_{l_2}} \EE F_{\pi_{l_3}} F_{\pi_{l_4}} \leq C n^{-q(\pi_{l_1}) +1 -q(\pi_{l_2},\pi_{l_3},\pi_{l_4}) + \epsilon(\pi)},\end{equation}
where $\epsilon(\pi) \in \{ 1, 2\}$, $\epsilon (\pi)= 2$ if there exists $j \in \{2,3,4\}$ such that $N(\pi_{l_j},\pi_{l_k}) = 0$ for all $k \in \{2,3,4\} \backslash \{j\}$, otherwise, $\epsilon(\pi) = 1$.

If $\epsilon(\pi)  = 1$ then  since there are at most $Cn^{q+q'}$ quadruples such that  $q(\pi_{l_1}) = q$ and $q(\pi_{l_2},\pi_{l_3},\pi_{l_4}) = q'$, we deduce  that $\sum_{(\pi_1,\cdots,\pi_4) \in S\cup S'} \ind(\epsilon(\pi)  =1  ) \EE F_{\pi_{l_1}}  F_{\pi_{l_2}} \EE F_{\pi_{l_3}} F_{\pi_{l_4}} \leq Cn^{2}$.

If $\epsilon(\pi)  = 2$, then, without loss of generality, we may assume  $N(\pi_{l_2},\pi_{l_k}) = 0$  for  $k \in \{3,4\}$. It implies that  $q(\pi_{l_2},\pi_{l_3},\pi_{l_4}) = q(\pi_{l_2}) + q(\pi_{l_3},\pi_{l_4})$. Since $(\pi_1,\cdots,\pi_4) \in S\cup S'$, $N(\pi_{l_2},\pi_{l_1}) \geq 1$, therefore $q(\pi_{l_1}) + q(\pi_{l_2}) \geq q(\pi_{l_1},\pi_{l_1}) +1$ and by Inequality (\ref{eq:type3}),
$$\EE F_{\pi_{l_1}} \EE  F_{\pi_{l_2}} \EE F_{\pi_{l_3}} F_{\pi_{l_4}} \leq C n^{-q(\pi_{l_1}) +1 -q(\pi_{l_2},\pi_{l_3},\pi_{l_4}) + 2} \leq C  n^{-q(\pi_{l_1},\pi_{l_2}) - q(\pi_{l_3},\pi_{l_4}) + 2}.$$
Finally, we notice that  there are at most $Cn^{q+q'}$ quadruples such that  $q(\pi_{l_1},\pi_{l_2}) = q$ and $q(\pi_{l_3},\pi_{l_4}) = q'$, it follows that:
$$\sum_{(\pi_1,\cdots,\pi_4)\in  S\cup S'} \EE F_{\pi_{l_1}}  \EE  F_{\pi_{l_2}}  F_{\pi_{l_3}} F_{\pi_{l_4}} \leq Cn^2.$$
Inequality (\ref{eq:4moment}) is proved.

\subsection{Proof of Theorem \ref{th:alphamain}}

By Equation (\ref{eq:momentnu}), for $m \geq 2$, we have (with $y_1 = 0$):
\begin{equation}
\label{eq:momentnu2} \nu_\gamma (P_m) =   f(0)^ m   +  \sum_{p =
2}^{m} \frac{ \gamma^{p-1}}{p!} \sum_{\phi \in \Sigma_{m,p}}
\int_{(\R^d)^{p-1}} \prod_{j=1}^ m f (y_{\phi(j)} - y_{\phi(j+1)})
dy_2 \cdots dy_{p}.
\end{equation}

The leading term in $\gamma$ is of order $\gamma^{m-1}$. Taking $p = m$ in the above expression gives:
$$
\nu_\gamma (P_m) \sim \gamma^{m-1} \int_{(\R^d)^{m-1}} \prod_{j=1}^ m f (y_{j} - y_{j+1}) dy_2 \cdots dy_{m}.
$$
A direct iteration leads to:
$$\int_{(\R^d)^{m-1}} \prod_{j=1}^ m f (y_{j} - y_{j+1}) dy_2 \cdots dy_{m} =  f^{*m} (0),$$
where $ f* g (y) = \int_{\R^d} f(x) g(y-x) dx$, $f^{*1} (x) = f(x)$, and for $m \geq 2$, $f^{*m} = f^{*(m-1)} * f$.

Hence $\int_{(\R^d)^{m-1}} \prod_{j=1}^ m f (y_{j} - y_{j+1}) dy_2 \cdots dy_{m} = \int_{\R^d} \hat f^m (\xi) d\xi = \int t^m \psi(t) dt$ and for all $m \geq 2$,
$$ \nu_\gamma (P_m) \sim \gamma^{m-1} \int  t^m \psi(t) dt = \int t^m \gamma ^{-2} \psi (\frac t \gamma) dt.$$
Since $\int t \psi(t) dt = \int \hat f(\xi) d\xi = f(0)$, this formula is still valid for $m =  1$. Now, let $h(t) = \sum_{ m \geq 1} h_m t^m$ with $\sum_{ m \geq 1} |h_m| t ^ m $ finite for all $t$, then since, $|\nu_\gamma ( P_m)| \leq m \gamma^ {m-1}  C^m$, using Fubini's Theorem, the conclusion follows.  \ep

\noindent {\bf Remark.}  
We can easily identify the next term in the asymptotic of $ \nu_\gamma (P_m)$. The second leading term in Equation (\ref{eq:momentnu2}) is of order $\gamma^{m-2}$. As in (\ref{eq:lemmap2}), it is equal to 
$$
I_m = \gamma^{m-2} \sum_{p=1} ^ {m-1} p \int_{\R} u ^ p \psi(u)
du   \int_{\R}  v ^ {m-p} \psi(v) dv.
$$
Since if $u \ne v$,  $\sum_{p=1} ^ {m-1}  p u ^p v^{m-p} = uv  (u - v )^{-2} ( (m-1) u^{m} - m u^{m-1} v + v ^m )$, we deduce that:
\begin{eqnarray*}
I_m & = & \gamma^{m-2} \int_{\R^2}   uv \frac{ (m-1) u^{m} - m u^{m-1} v + v ^m } {  (u - v)^2 }  \psi(u) \psi(v) du  dv.
\end{eqnarray*}

\subsection{Proof of Proposition \ref{prop:upperbound}}

Let $D_n$ denote the $n \times n$ matrix with entry $i,j$ equal to: $\ind( \| X_ i - X_j \| \leq \delta_n)$, if $I_n$ denotes the $n \times n$ identity matrix, $D_n - I_n$ is the adjacency matrix of the random geometric graph $\mathcal G ( {\cX_n }, \delta_n)$ where there is an edge between $i \ne j$ if $\| X_ i - X_j \| \leq \delta_n$. We have component wise:
$$
- \| f \|_{\infty} D_n \leq B_n \leq \| f \|_{\infty} D_n.
$$
Since the spectral radius $\rho(B_n)$  of $B_n$ is upper bounded by $\max_{ 1 \leq i \leq n} | \sum_{j=1 }^n (B_n )_{ij} |$, we deduce that:
$$
\rho(B_n) \leq \| f \|_{\infty} (1+ \Delta_n),
$$ 
where $\Delta_n$ is the maximal degree of the graph $\mathcal G( {\cX_n }, \delta_n)$. Then, the proposition follows from Theorem 6.6 of Penrose \cite{penrose}. \ep

\section{Further properties of the Euclidean Random Matrices}
\label{sec:misc}

\subsection{Eigenvectors of Euclidean Random Matrices} 

As it is  pointed by  M\'ezard, Parisi and Zee \cite{mezardparisi}, if $U_i = (\Phi_{k,n})_i = e^{2i \pi k.X_i}$ we have:
\begin{eqnarray}
\label{eq:pseudoeig}
(A \Phi_{k,n})_i & = & \Bigm( \sum_j F(X_i - X_j) e^{-2 i \pi k. (X_i - X_j )  } \Bigm) (\Phi_{k,n})_i,
\end{eqnarray}

In particular, if $F(x) = e^{2 i \pi k. x}$, then  $n$ is an eigenvalue with $\Phi_{k,n}$ as eigenvector and the rank of $A$ is $1$. Note also by the Strong Law of Large Numbers that for all $i$, a.s.
$$
\lim_{ n \to \infty} \frac{1}{ n} \sum_j F(X_i - X_j) e^{-2 i \pi k. (X_i - X_j)} =  \hat F(k).
$$
if $A_n =  A / n$, by Equation (\ref{eq:pseudoeig}), for all $i$, a.s.:
$$
\lim_{n \to \infty}  (A_n \Phi_{k,n})_i = \hat F(k) (\Phi_{k,n})_i.
$$
This last equation is consistent with Theorem \ref{th:main}:  a.s.  for $n$ large enough there exists an eigenvalue of $A_n$ close to $\hat F(k)$. It is possible to strengthen this last convergence as follows:
\begin{proposition}
For $p \geq 1$, let $\| U \|_p = \bigm(\sum_{i\geq 1} |U_i|^p \bigm)^{1/p}$ and $\| U \|_\infty = \sup_{i\geq 1} |U_i|$. For all $p \in ( 2, \infty]$,  a.s. for all $k \in \Z^d$,
$$
\lim_{n \to \infty} \| A_n \Phi_{k,n} -\hat F(k) \Phi_{k,n} \|_p = 0.
$$
Moreover,
$\lim_{n \to \infty} \EE \| A_n \Phi_{k,n} -\hat F(k) \Phi_{k,n} \|^2_2 = \| f \|^2 _2 - |\hat F(k)|^ 2 = \sum_{l \neq k}  |\hat F(l)|^ 2$.
\end{proposition}

\bp
To simplify notation, we write $\Phi = \Phi_{k,n}$ and $f_{ij} = F(X_i - X_j)$.
\begin{eqnarray*}
\PP \Bigm( \| A_n \Phi_{k,n} -\hat F(k) \Phi_{k,n} \|_p > \epsilon \Bigm) & = & \PP \Bigm( \sum_{i =1 } ^ n \Bigm | \frac 1 n \sum_{j=1} ^ n  f_{ij} \Phi_j - \hat F (k) \Phi_i \Bigm |^p > \epsilon ^p  \Bigm) \\
& \leq & n \PP \Bigm(  \Bigm | \frac 1 n \sum_{j=1} ^ nf_{1j} \Phi_j - \hat F (k) \Phi_1 \Bigm |^p  >  \frac {\epsilon^p}{n}\Bigm) \\
& \leq & n \PP \Bigm(  \Bigm |\sum_{j=1} ^ n f_{1j} \Phi_j - n \hat F (k) \Phi_1 \Bigm |  >  \epsilon  n^{1- 1/p}\Bigm).
\end{eqnarray*}
From Equation (\ref{eq:pseudoeig}), $|\sum_{j=1} ^ n f_{1j} \Phi_j - n \hat F (k) \Phi_1 | =  |\sum_{j=1} ^ n F(X_1 - X_j) e^{-2 i \pi k. (X_1 - X_j )}  - n \hat F (k)|$. Hence:
\begin{eqnarray*}
\PP \Bigm( \| A_n \Phi_{k,n} -\hat F(k) \Phi_{k,n} \|_p > \epsilon \Bigm) & \leq & n  \PP \Bigm( | \sum_{j=1} ^ n F(X_1 - X_j) e^{-2 i \pi k. (X_1 - X_j )} - n \hat F(k) |  >  \epsilon  n^{1- 1/p} \Bigm ) \\
 & \leq & n \EE\Bigm [ \PP \Bigm ( | \sum_{j=2}^ n  F(X_1 - X_j) e^{-2 i \pi k. (X_1 - X_j )} - (n-1) \hat F(k) |   \\
&& \quad \quad \quad\quad  \quad \quad\quad  \quad >  \epsilon  n^{1- 1/p} - |F(0)| - |\hat F(k)|  \Bigm |  X_1\Bigm )\Bigm ] \\
& \leq & 2 n \exp ( - \frac { \max(0,( \epsilon n^{1-1/p} - |F(0)| - |\hat F (k)|  ))^2 } { \| F \|_{\infty} (n-1)}),
\end{eqnarray*}
where the last equation is Hoeffding's Inequality. We then apply Borel Cantelli Lemma.

It remains to prove the statement of the proposition for $p = 2$. Similarly, we obtain:
$$
\EE \| A_n \Phi_{k,n} -\hat F(k) \Phi_{k,n} \|^2_2 = \frac {1}{n} \EE |\sum_j F(X_1 - X_j) e^{-2 i \pi k. (X_1 - X_j )}  - n \hat F (k)|^2. $$
We then write $\EE |\sum_j F(X_1 - X_j) e^{-2 i \pi k. (X_1 - X_j )}  - n \hat F (k)|^2 = \EE [ \EE [ |\sum_j \bigm ( F(X_1 - X_j) e^{-2 i \pi k. (X_1 - X_j )}  -  \hat F (k) \bigm) |^2 | X_1 ]] = |F(0) - \hat F (k)| ^ 2 + (n-1) \int_\Omega | F(x) e^{-2i \pi k.x} - \hat F(k)|^2 dx$. The statement follows.
\ep

\subsection{Correlation of the Eigenvalues}

In this paragraph, we state an elementary lemma on the $m$-correlation of the eigenvalues of $A$ ($ m \leq n)$:
$$
M_m = 1/ \binom{n}{m}  \EE \sum_{\{i_1,\cdots,i_m\}\subset \{1,\cdots n\}} \prod_{j=1}^m (\lambda_{i_j} - F(0)),
$$
where the sum is over all subsets of $\{1,\cdots,n\}$ of cardinal $m$. Note that $M_1 = 0$ and that $M_m$ is related to the factorial moment measure $\rho_m(dz_1,\cdots,dz_m)$ (also called the joint intensity measure, refer to Daley and Vere-Jones \cite{daley}) of the point process $\{\lambda_{1} - F(0),\cdots, \lambda_n - F(0) \}$  as follows:
$$
M_m = \int_{\C^m} \prod_{j=1} ^ m z_j \rho_m (dz_1,\cdots,dz_m),
$$
Heuristically, $\rho_m(dz_1,\cdots,dz_m)$ is the infinitesimal probability of having an eigenvalue at $F(0) + z_i$ for each $i \in \{1,\cdots, m\}$. We define $\bar A = A - F(0) I$, where $I$ is the $n \times n$ identity matrix (note that ``$\bar \cdot$`` is not the complex conjugate of the matrix $A$). $\bar A(x_1,\cdots,x_m)$ is the $m \times m$ matrix where the coefficient $i,j$ is equal to $F(x_i-x_j) - \delta_{ij} F(0)$.

\begin{lemma}
\label{eq:Cm}
$$
M_{m} =  \int_{\Omega^{m}} \det \bar A (x_1,\cdots,x_{m}) dx_1 \cdots dx_{m}.
$$
\end{lemma}
 For $m = 2$ we get:
$$
M_2 = - \int_{\Omega} F(x)^2 dx,
$$
the point process of eigenvalues is thus repulsive.

\bp 
The characteristic polynomial of $\bar A$ is $\chi_{\bar A} (t) = \det (\bar A-tI) = \prod_{i=1} ^ n (\lambda_{i} - f(0) - t) = \sum_{m=0} ^ n a_m  (-t) ^{n-m}$, where, $a_m = \sum_{\{i_1,\cdots,i_m\}} \prod_{j=1}^m (\lambda_{i_j} - f(0))$.  However, by Newton formula, we also have, $a_m = \sum_{\{i_1,\cdots,i_m\}} det \bar A_{\{i_1,\cdots,i_m\}}$, where for a set of indices ${\bf i} = \{i_1,\cdots,i_m\}$, $\bar A_{\bf i}$ is the $m \times m$ extracted matrix obtained from $\bar A$ by keeping the raws and columns $\{i_1,\cdots,i_m\}$ (i.e. $\bar A_{\bf i}$ is a principal minor). Taking expectation, we deduce that $\EE a_m = \binom{n}{m} \int_{\Omega^m} \det \bar A (x_1,\cdots,x_{m}) dx_1 \cdots dx_{m}$.\ep

In Lemma \ref{eq:Cm}, we have computed the mean value of the symmetric polynomials: $$\alpha_m (x_1,\cdots,x_n) = \sum_{\{i_1,\cdots,i_m\}\subset \{1,\cdots n\}} \prod_{j=1}^m x_j$$
for the vector ${\bf \bar \lambda} = (\lambda_{1} - F(0),\cdots, \lambda_n - F(0)).$ Actually, it is possible to compute the mean value of the symmetric polynomials:
$$
\alpha_{m,k} (x_1,\cdots,x_n) =  \sum_{\{i_1,\cdots,i_m\} \subset \{1,\cdots n\}} \prod_{j=1}^m x^k_j.
$$
for the vector ${\bf \bar \lambda}$.  To this end simply consider, $\chi_{\bar A^k} (t) = \det (\bar A^k -tI) = \prod_{i=1} ^ n ((\lambda_{i} - f(0))^k - t)$. We obtain similarly:
$$
 1/ \binom{n}{m} \EE \alpha_{m,k} (\bar \lambda) = \int_{\C^n} \prod_{j=1} ^ m z_j ^ k \rho_m (dz_1,\cdots,dz_m) = \int_{\Omega^n} \det \bar A^k_{{\bf m }} (x_1, \cdots, x_n) dx_1 \cdots dx_n,
$$
where ${\bf  m } = \{1, \cdots,m\}$ and  for a set of indices ${\bf i} = \{i_1,\cdots,i_m\}$, $\bar A_{\bf i}$ is the $m \times m$ extracted matrix obtained from $\bar A$ by keeping the raws and columns $\{i_1,\cdots,i_m\}$.

\subsection*{Acknowledgments}

The author thanks Neil O'Connell for suggesting this problem and Florent Benaych-George for valuable comments. 

\bibliographystyle{plain}

\bibliography{../../bib}

\begin{thebibliography}{10}

\bibitem{bai}
Z.~D. Bai.
\newblock Methodologies in spectral analysis of large-dimensional random
  matrices, a review.
\newblock {\em Statist. Sinica}, 9(3):611--677, 1999.
\newblock With comments by G. J.\ Rodgers and Jack W.\ Silverstein; and a
  rejoinder by the author.

\bibitem{bogomolny}
E.~Bogomolny, O.~Bohigas, and C.~Schmidt.
\newblock Spectral properties of distance matrices.
\newblock {\em Journal of Physics A: Mathematical and General}, 36:3595--3616,
  2003.

\bibitem{dembo2006}
W.~Bryc, A.~Dembo, and T.~Jiang.
\newblock Spectral measure of large random {H}ankel, {M}arkov and {T}oeplitz
  matrices.
\newblock {\em Ann. Probab.}, 34(1):1--38, 2006.

\bibitem{fanchung03}
F.~Chung, L.~Lu, and V.~Vu.
\newblock Eigenvalues of random power law graphs.
\newblock {\em Ann. Comb.}, 7(1):21--33, 2003.

\bibitem{fanchung04}
F.~Chung, L.~Lu, and V.~Vu.
\newblock The spectra of random graphs with given expected degrees.
\newblock {\em Internet Math.}, 1(3):257--275, 2004.

\bibitem{daley}
D.J. Daley and D.~Vere-Jones.
\newblock {\em An introduction to the Theory of Point Processes}.
\newblock Springer Series in Statistics. Springer-Verlag, New-York, 1988.

\bibitem{ganesh05}
M.~Draief and A.~Ganesh.
\newblock Efficient routeing in {P}oisson small-world networks.
\newblock {\em J. Appl. Probab.}, 43(3):678--686, 2006.

\bibitem{draief}
M.~Draief, A.~Ganseh, and L.~Massouli\'e.
\newblock Thresholds for virus spread on networks.
\newblock {\em to appear in Ann. Appl. Probab.}

\bibitem{furedikomlos}
Z.~F{\"u}redi and J.~Koml{\'o}s.
\newblock The eigenvalues of random symmetric matrices.
\newblock {\em Combinatorica}, 1(3):233--241, 1981.

\bibitem{hammond}
C.~Hammond and S.~Miller.
\newblock Distribution of eigenvalues for the ensemble of real symmetric
  {T}oeplitz matrices.
\newblock {\em J. Theoret. Probab.}, 18(3):537--566, 2005.

\bibitem{gine}
V.~Koltchinskii and E.~Gin{\'e}.
\newblock Random matrix approximation of spectra of integral operators.
\newblock {\em Bernoulli}, 6(1):113--167, 2000.

\bibitem{Lovasz}
L.~Lov{\'a}sz.
\newblock Random walks on graphs: a survey.
\newblock In {\em Combinatorics, Paul Erd\H os is eighty, Vol.\ 2 (Keszthely,
  1993)}, volume~2 of {\em Bolyai Soc. Math. Stud.}, pages 353--397. J\'anos
  Bolyai Math. Soc., Budapest, 1996.

\bibitem{mezardparisi}
M.~M{\'e}zard, G.~Parisi, and A.~Zee.
\newblock Spectra of {E}uclidean random matrices.
\newblock {\em Nuclear Phys. B}, 559(3):689--701, 1999.

\bibitem{offersimons}
C.~Offer and B.~D. Simons.
\newblock Field theory of {E}uclidean matrix ensembles.
\newblock {\em J. Phys. A}, 33(42):7567--7583, 2000.

\bibitem{penrose}
M.~Penrose.
\newblock {\em Random Geometric Graphs}.
\newblock Oxford Studies in Probability. Oxford Univeristy Press, Oxford, 2003.

\bibitem{talagrand96}
M.~Talagrand.
\newblock A new look at independence.
\newblock {\em Ann. Probab.}, 24(1):1--34, 1996.

\bibitem{vershik}
A.~M. Vershik.
\newblock Random metric spaces and universality.
\newblock {\em Russian Math. Surveys}, 59(2):259--295, 2004.

\bibitem{wigner}
E.~Wigner.
\newblock On the distribution of the roots of certain symmetric matrices.
\newblock {\em Ann. of Math. (2)}, 67:325--327, 1958.

\end{thebibliography}

\end{document}